\documentclass[runningheads]{llncs}

\usepackage{amssymb}
\usepackage{amsmath}
\usepackage{graphicx}
\usepackage{xcolor}
\usepackage{hyperref}
\usepackage{caption}
\usepackage{listings}
\usepackage{booktabs}
\usepackage{array}

\usepackage{tabularx}

\usepackage{makecell}

\usepackage{multirow}

\urldef{\mailsb}\path|jan.valdman@utia.cas.cz |

\definecolor{myart}{RGB}{0,0,0}

\def\BibTeX{{\rm B\kern-.05em{\sc i\kern-.025em b}\kern-.08em
    T\kern-.1667em\lower.7ex\hbox{E}\kern-.125emX}}
    \usepackage{comment}

\usepackage{ulem}

\newcommand\dx{\, \mathrm{d}x}
\newcommand\dt{\, \mathrm{d}t}
\newcommand\ds{\, \mathrm{d}s}

\newcommand{\green}{\color{teal}}

\newcommand{\uncom}[1]{{\green #1}}

\definecolor{defretcolor}{rgb}{0.00, 0.42, 0.24} 
\definecolor{jnpcolor}{rgb}{0.86, 0.08, 0.24} 
\definecolor{funccolor}{rgb}{0.06, 0.25, 0.49} 
\definecolor{numbercolor}{rgb}{0.48, 0.00, 0.48} 
\definecolor{commentcolor}{rgb}{0.38, 0.63, 0.69} 
\definecolor{backcolour}{rgb}{0.95, 0.95, 0.96} 
\definecolor{stringcolor}{rgb}{0.31, 0.60, 0.02} 

\lstdefinestyle{mystyle}{
    backgroundcolor=\color{backcolour},
    commentstyle=\color{commentcolor},
    basicstyle=\ttfamily\small,
    breakatwhitespace=false,
    breaklines=true,
    captionpos=b,
    keepspaces=true,
    numbers=left,
    numbersep=5pt,
    showspaces=false,
    showstringspaces=false,
    showtabs=false,
    tabsize=2,
    keywordstyle=\color{black}, 
    stringstyle=\color{stringcolor},
    classoffset=0, 
    morekeywords={def, return}, 
    keywordstyle=\color{defretcolor},
    emph={array, sum, at, abs, log, set, dot}, 
    emphstyle=\color{funccolor},
    classoffset=0, 
    morekeywords={jnp},
    keywordstyle=\color{jnpcolor},
    classoffset=2, 
    morekeywords={0,1,2,3,4,5,6,7,8,9},
    keywordstyle=\color{numbercolor}, 
    classoffset=0, 
}

\title{Physics-Informed Neural Network for Diffusion-Reaction Problems with Dead-Core Formation in Catalyst Slabs
}

 \titlerunning{Physics-Informed Neural Network for Diffusion-Reaction Problems}

\author{Piotr Skrzypacz\inst{1}\orcidID{0000-0001-6422-5469}, 
Kaisar Tangirbergen\inst{1},
Jan Valdman\inst{2}\orcidID{0000-0002-6081-5362}}
\authorrunning{Piotr Skrzypacz, Kaisar Tangirbergen, Jan Valdman}
\institute{School of Sciences and Humanities, Nazarbayev University, 53 Kabanbay Batyr Ave., Astana 010000, Kazakhstan \\
\email{piotr.skrzypacz@nu.edu.kz; kaisar.tangirbergen@nu.edu.kz}
\and
Department of Mathematics, Faculty of Science, University of South Bohemia, 
Brani\v sovsk\' a 31, 37005~\v{C}esk\'{e}~Bud\v{e}jovice, Czech Republic
\\
\email{jvaldman@prf.jcu.cz}
}
\date{}

\begin{document}

\maketitle

\begin{abstract}
This work investigates a nonlinear two-point boundary value problem
arising in diffusion--reaction processes in catalyst slabs with
power-law kinetics and fractional reaction order. For sufficiently
large Thiele modulus, the solution develops a dead-core region separated
from the active region by an unknown free boundary. We propose a
structured Physics-Informed Neural Network (PINN) framework that
incorporates the asymptotic behavior at the dead-core interface into a
hard-constrained trial solution and treats the interface location as a
trainable parameter. The concentration profile and free boundary are
therefore approximated simultaneously without explicit interface
tracking or penalty-based enforcement of the interface conditions. The
method is validated against the exact solution for power-law kinetics
and a high-precision numerical shooting method. Numerical experiments
covering near-critical and strongly supercritical regimes demonstrate
accurate recovery of both the concentration profile and dead-core
location, together with robustness to random initialization and
collocation sampling. While classical shooting is more efficient for
the present one-dimensional benchmark, the proposed formulation provides
a flexible framework for extensions to multidimensional geometries and
problems for which analytical solutions are unavailable.
\end{abstract}

\section{Introduction}

Nonlinear diffusion--reaction equations arise in a variety of physical
and engineering applications, including transport and reaction in porous
catalysts and biological systems \cite{Non,rec-diff}. In catalytic
materials, the interplay between diffusion and reaction can produce
regions in which the reactant concentration vanishes. For fractional-order
kinetics, such regions are commonly referred to as dead cores or dead
zones \cite{Aris,Bandle,Temkin}. Dead-core formation has been considered
in a range of catalytic and biochemical processes, and it can substantially
affect concentration profiles and performance measures such as the
effectiveness factor of a catalytic particle \cite{york2011dead}.

Dead-core diffusion--reaction problems are challenging numerically because
the fractional reaction term is non-Lipschitz at vanishing concentration.
This behavior can make conventional numerical treatment more delicate
near the free boundary \cite{aziz1988numerical,barrett1991finite}.
Various numerical and semi-analytical approaches have been developed for
dead-core problems with power-law kinetics and related catalyst models
\cite{Dead-core,membrane_piotr,Proff2023}. For catalyst slabs, the
critical Thiele modulus provides a natural criterion for the onset of
dead-core formation and relates the diffusion and reaction rates to the
emergence of the free boundary \cite{Prof2022,Prof2022_general}.

In this work, we consider the power-law kinetics
\begin{equation}\label{LH_kin}
R(C_A)=k C_A^n,
\end{equation}
where $C_A$ is the molar concentration of species $A$, $k$ is the reaction
rate constant, and $n\in[0,1)$ is the reaction exponent. We use the
convention
\[
C_A^n:=
\begin{cases}
\bigl(\max\{C_A,0\}\bigr)^n,
& n\in(0,1),\\[1.2ex]
\displaystyle\operatorname{sign}\bigl(\max\{C_A,0\}\bigr),
& n=0,
\end{cases}
\]
where $\operatorname{sign}(t)=t/|t|$ for $t\neq0$ and
$\operatorname{sign}(0)=0$.

The main contribution of this work is a structured free-boundary
Physics-Informed Neural Network (PINN) formulation for dead-core
diffusion--reaction problems. The approach combines the analytically
derived local asymptotic behavior at the dead-core interface with a
hard-constrained trial solution, while the unknown interface location is
represented by a trainable scalar parameter. The concentration profile
and free boundary are therefore identified simultaneously, with the
interface and outer boundary conditions satisfied by construction rather
than through penalty terms.

The use of a PINN is motivated here not by computational superiority over
classical one-dimensional solvers, but by the flexibility of the
formulation. The present slab problem admits an exact solution and is
therefore particularly suitable as a controlled benchmark: it allows the
accuracy of the learned concentration profile and free-boundary location
to be assessed independently of the neural-network training procedure.
At the same time, the structured representation provides a formulation
that can be extended to settings in which explicit analytical solutions
are unavailable or classical one-dimensional shooting becomes less
applicable. In particular, the approach is intended as a basis for
multidimensional catalyst geometries, inverse identification of kinetic
and transport parameters, and parametric nonlinear diffusion--reaction
problems with unknown interfaces.

The remainder of the paper first introduces the mathematical model and
analytical dead-core structure that motivate the PINN construction. The
free-boundary PINN and its training procedure are then described, followed
by numerical validation, comparison with a classical shooting method,
robustness tests, and ablation experiments.

\section{Mathematical model for dead-core formation in catalyst slabs}

\subsection{Mass balance equation in catalyst slabs}

We consider a catalyst slab of half-length $R_p$ containing a single
reaction. Let $r_p$ denote the distance from the slab center. Assuming
standard Fickian diffusion, the flux is
\begin{equation}\label{eq_diff_flux}
j(r_p)=-D_{\mathrm{eff}}\frac{dC_A(r_p)}{dr_p},
\end{equation}
where $C_A$ is the molar concentration of species $A$ and
$D_{\mathrm{eff}}$ is the effective diffusion coefficient. The steady
diffusion--reaction equation is
\begin{equation}\label{eq_diff_react}
D_{\mathrm{eff}}\frac{d^2C_A}{dr_p^2}
=
R(C_A),
\qquad r_p\in(0,R_p),
\end{equation}
where $R(C_A)$ is given by the power-law kinetics~\eqref{LH_kin}.
The symmetry and bulk-concentration conditions are
\begin{equation}\label{eq_bc_model}
\frac{dC_A}{dr_p}(0)=0,
\qquad
C_A(R_p)=C_{A,b},
\end{equation}
with $C_{A,b}>0$. Integrating~\eqref{eq_diff_react} from $0$ to $r_p$
gives
\begin{equation}\label{eq_dCdr}
D_{\mathrm{eff}}\frac{dC_A}{dr_p}(r_p)
=
\int_0^{r_p} kC_A^n(t)\dt\geq0.
\end{equation}
Thus $C_A$ is monotonically non-decreasing, with zero derivative at the
slab center and $C_A(R_p)=C_{A,b}$.

\subsection{Dimensionless problem}

Introducing the dimensionless distance and concentration
\[
x=\frac{r_p}{R_p},
\qquad
u=\frac{C_A}{C_{A,b}},
\]
and the Thiele modulus
\begin{equation}\label{eq_Thiele}
\phi=
\sqrt{\frac{R_p^2 C_{A,b}^{\,n-1}k}{D_{\mathrm{eff}}}},
\end{equation}
the problem becomes
\begin{equation}\label{eq_nondim}
u''=\phi^2r_n(u),
\end{equation}
where
\begin{equation}\label{eq_ru}
r_n(u)=u^n.
\end{equation}
The boundary conditions are
\begin{equation}\label{eq_bc_nondim}
u'(0)=0,
\qquad
u(1)=1.
\end{equation}

\subsection{Dead-core length}

We first consider the more general problem
\begin{equation}\label{eq_twopoint_general}
\begin{split}
u''&=\phi^2r(u)\quad\text{in }(0,1),\\
u'(0)&=0,\qquad u(1)=1,
\end{split}
\end{equation}
where $r:[0,\infty)\to[0,\infty)$ is nonnegative and continuously differentiable on $(0,\infty)$, with $r(u)>0$ for $u>0$. Suppose that for $\phi>\phi^*$ the
solution develops a nontrivial dead core,
\[
u(x)=0\quad\text{for }0\leq x\leq x_{dz},
\qquad
u(x)>0\quad\text{for }x_{dz}<x\leq1.
\]
At the interface,
\[
u(x_{dz})=u'(x_{dz})=0.
\]
Since $u''\geq0$, the solution is convex and the dead core is an interval
adjacent to the symmetry point.

On the active region, multiplying~\eqref{eq_twopoint_general} by $u'$
and integrating from $x_{dz}$ to $x$ yields the first integral
\begin{equation}\label{eq_first_integral}
\frac{1}{2}\bigl(u'(x)\bigr)^2
=
\phi^2\int_0^{u(x)}r(t)\dt.
\end{equation}
The integration constant vanishes because both $u$ and $u'$ are zero at
the dead-core interface. Since $u$ is increasing on the active region,
Eq.~\eqref{eq_first_integral} gives
\[
\frac{dx}{du}
=
\frac{1}{\sqrt{2}\phi}
\left(
\int_0^u r(t)\dt
\right)^{-1/2}.
\]
Integrating from $u=0$ at $x=x_{dz}$ to $u=1$ at $x=1$ gives
\begin{equation}\label{eq_active_length_general}
1-x_{dz}
=
\frac{1}{\sqrt{2}\phi}
\int_0^1
\left(
\int_0^s r(t)\dt
\right)^{-1/2}\ds.
\end{equation}
This motivates the critical Thiele modulus
\begin{equation}\label{eq_general_phistar}
\phi^*
=
\frac{1}{\sqrt{2}}
\int_0^1
\left(
\int_0^s r(t)\dt
\right)^{-1/2}\ds,
\end{equation}
provided the integral is finite. Consequently,
\begin{equation}\label{eq_xdz_length}
x_{dz}=1-\frac{\phi^*}{\phi}.
\end{equation}
Thus, for fixed kinetics, the length of the active region is
\[
1-x_{dz}=\frac{\phi^*}{\phi},
\]
so it decreases inversely with the Thiele modulus.

\begin{figure}[t]
    \centering
    \includegraphics[width=1\linewidth]{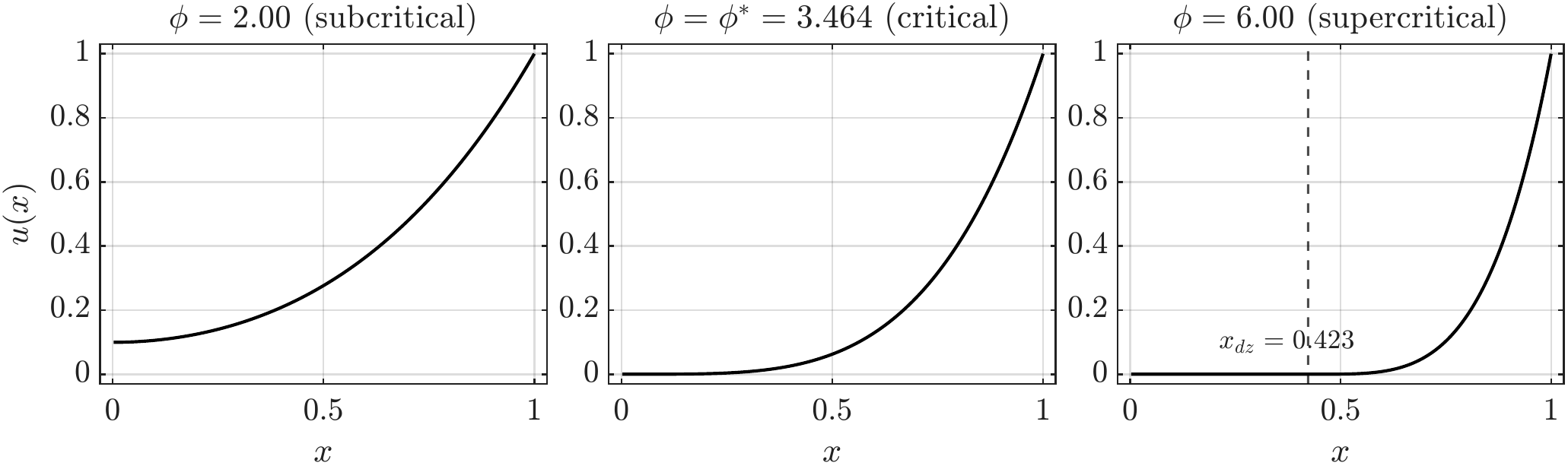}
    \caption{Exact solutions $u(x)$ for $n=0.5$ and
    $\phi\in\{2,2\sqrt{3},6\}$. The critical value is
    $\phi^*=2\sqrt{3}\approx3.464$, and for $\phi=6$ the dead-core
    boundary is $x_{dz}=1-1/\sqrt{3}\approx0.42265$.}
    \label{fig:solution_example}
\end{figure}

For the slab geometry considered here, the solution is illustrated in
Fig.~\ref{fig:solution_example}. Below the critical value $\phi^*$ the
solution is positive throughout the slab, whereas for $\phi>\phi^*$ a
finite dead core develops.

For the power-law kinetics~\eqref{eq_ru},
\[
\int_0^s t^n\dt=\frac{s^{n+1}}{n+1},
\]
and Eq.~\eqref{eq_general_phistar} gives the explicit critical Thiele
modulus
\begin{equation}\label{formula_phistar}
\phi^*=\frac{\sqrt{2(1+n)}}{1-n}.
\end{equation}
Substitution into~\eqref{eq_xdz_length} yields
\begin{equation}\label{eq_xdz_power}
x_{dz}
=
1-\frac{\sqrt{2(1+n)}}{\phi(1-n)}.
\end{equation}

The corresponding exact solution can be obtained directly from the first
integral~\eqref{eq_first_integral}. It is
\begin{equation}\label{eq_udeadcore_power}
u(x)=
\begin{cases}
0,
&x\in[0,x_{dz}],\\[0.8ex]
\displaystyle
\left(
\frac{x-x_{dz}}{1-x_{dz}}
\right)^{\frac{2}{1-n}},
&x\in[x_{dz},1].
\end{cases}
\end{equation}
At $\phi=\phi^*$, Eq.~\eqref{eq_xdz_power} gives $x_{dz}=0$ and
Eq.~\eqref{eq_udeadcore_power} reduces to the critical solution. For
$\phi>\phi^*$, the solution contains a nontrivial dead-core region.
Equations~\eqref{eq_xdz_power} and~\eqref{eq_udeadcore_power} provide
the analytical reference for validating the proposed PINN.

\section{Physics-Informed Neural Network Approach}

To solve the nonlinear diffusion--reaction problem with an unknown
dead-core interface, we employ a structured free-boundary
Physics-Informed Neural Network (PINN). The objective is to approximate
the concentration profile in the active region while simultaneously
identifying the unknown dead-core location. The analytical structure
derived in Section~2 is incorporated directly into the neural-network
trial solution.

\subsection{Domain Transformation}

For $\phi>\phi^*$, the solution is zero on $[0,x_{dz}]$ and positive on
$(x_{dz},1]$. We therefore restrict the neural-network approximation to
the active region and introduce the transformed coordinate
\begin{equation}\label{eq_xi_transform}
\xi=\frac{x-x_{dz}}{1-x_{dz}},
\qquad
x\in[x_{dz},1],
\qquad
\xi\in[0,1].
\end{equation}
Thus, the unknown dead-core interface is mapped to $\xi=0$ and the slab
surface to $\xi=1$.

The behavior of the solution near the dead-core interface provides
additional information that can be incorporated into the network
architecture. For $0<n<1$, assume locally that
\[
u(x)\sim C(x-x_{dz})^p,
\qquad C>0,
\qquad x\to x_{dz}^{+}.
\]
Substitution into
\[
u''=\phi^2u^n
\]
gives, to leading order,
\[
Cp(p-1)(x-x_{dz})^{p-2}
\sim
\phi^2C^n(x-x_{dz})^{pn}.
\]
Balancing the powers of $(x-x_{dz})$ yields
\begin{equation}\label{eq_asymptotic_p}
p-2=pn,
\qquad
p=\frac{2}{1-n}.
\end{equation}
For $0<n<1$, one has $p>2$, while for $n=0$, $p=2$; in both cases the resulting power satisfies
\[
u(x_{dz}) = u'(x_{dz}) = 0.
\]

We incorporate this leading-order behavior into the neural-network
trial solution by defining
\begin{equation}\label{eq_structured_ansatz}
u(\xi)
=
\xi^p
\left[
1+(1-\xi)\mathcal{M}_{\theta}(\xi)
\right]^2,
\qquad
p=\frac{2}{1-n},
\end{equation}
where $\mathcal{M}_{\theta}$ denotes the output of a feedforward neural
network with trainable parameters $\theta$.

The structured form in Eq.~\eqref{eq_structured_ansatz} imposes the
relevant boundary conditions without adding corresponding penalty terms
to the loss. At $\xi=0$, the factor $\xi^p$ gives
\[
u(0)=u'(0)=0,
\]
while at $\xi=1$ the correction term vanishes and therefore
\[
u(1)=1.
\]
Consequently, the network is required to learn the differential
equation in the active region rather than simultaneously learning the
boundary conditions.

\subsection{Neural-Network Architecture}

The correction function $\mathcal{M}_{\theta}$ is represented by a fully
connected feedforward neural network. The same architecture is used for
all test cases to avoid case-specific tuning. The network consists of
4 hidden layers with 64 neurons
per hidden layer, using hyperbolic tangent activation
functions, followed by a linear output layer.

The architecture is kept fixed while the reaction order and Thiele
modulus are varied. This allows the numerical experiments to assess the
effect of the physical parameters rather than changes in network
complexity. The trainable scalar parameter describing the free boundary,
introduced below, is optimized together with the network weights.

\subsection{Free-Boundary Parameterization}

The dead-core interface $x_{dz}$ is treated as an additional trainable
parameter. To ensure that it remains in the admissible interval, we
introduce an unconstrained scalar parameter $\alpha$ and define
\begin{equation}\label{eq_xdz_parameter}
x_{dz}
=
(1-\varepsilon)\sigma(\alpha),
\qquad
\sigma(\alpha)=\frac{1}{1+e^{-\alpha}},
\end{equation}
where $\varepsilon>0$ is a small numerical constant. This
parameterization guarantees
\[
0<x_{dz}<1-\varepsilon
\]
during training and prevents the active interval from collapsing to zero
length.

The interface location therefore enters the model in two ways: through
the transformed coordinate~\eqref{eq_xi_transform} and through the
scaling of the differential operator. Gradients of the loss with respect
to $\alpha$ consequently provide a direct mechanism for optimizing the
free-boundary location.

\subsection{Physics-Informed Loss Function}

Using the transformation~\eqref{eq_xi_transform}, the derivatives with
respect to $x$ and $\xi$ are related by
\[
\frac{d}{dx}
=
\frac{1}{1-x_{dz}}\frac{d}{d\xi},
\]
and hence
\begin{equation}\label{eq_second_transform}
\frac{d^2u}{dx^2}
=
\frac{1}{(1-x_{dz})^2}
\frac{d^2u}{d\xi^2}.
\end{equation}
The residual of the governing equation in the transformed active region
is therefore
\begin{equation}\label{eq_pinn_residual}
\mathcal{R}(\xi)
=
\frac{1}{(1-x_{dz})^2}
\frac{d^2u}{d\xi^2}
-
\phi^2u^n.
\end{equation}

The network parameters $\theta$ and the free-boundary parameter
$\alpha$ are obtained by minimizing the mean-squared residual
\begin{equation}\label{eq_pinn_loss}
\mathcal{L}_{\mathrm{PDE}}
=
\frac{1}{N}
\sum_{i=1}^{N}
\mathcal{R}(\xi_i)^2,
\end{equation}
where $\{\xi_i\}_{i=1}^{N}$ are collocation points in $[0,1]$.

Because the structured ansatz satisfies the interface and outer
boundary conditions by construction, no additional boundary-penalty
terms are required. The resulting optimization problem therefore
contains the neural-network parameters and the single free-boundary
parameter as the unknown quantities.

For numerical stability, the nonlinear term is evaluated using
\[
u_{\mathrm{safe}}=\max(u,10^{-30})
\]
when computing $u^n$. This safeguard is relevant for $0<n<1$, since
the derivative $nu^{n-1}$ becomes singular as $u\to0^+$. It is used only in the numerical evaluation of the residual and does not alter the analytical formulation.

\subsection{Sampling Strategy}

Collocation points are sampled in the transformed domain
$\xi\in[0,1]$. To resolve both the dead-core interface and the outer
boundary, we use a three-component sampling strategy. One third of the
points are sampled uniformly, while the remaining points are
concentrated near $\xi=0$ and $\xi=1$. Specifically,
\begin{equation}\label{eq_biased_sampling}
\xi_0=U^\beta,
\qquad
\xi_1=1-U^\beta,
\qquad
U\sim U[0,1],
\qquad
\beta>1.
\end{equation}
The uniform component provides global coverage, whereas the two biased
components increase the resolution near the regions where the solution
may vary most rapidly.

The same sampling strategy is used for all test cases. During the
warm-up and Adam stages, $N=3000$ collocation points are used with
$\beta=5$. During the L-BFGS refinement stages, fixed sets containing
$N=4000$ and $N=8000$ points are used, with $\beta=6$. The use of the
same sampling parameters for all cases avoids case-specific adjustment
of the numerical procedure.

\subsection{Optimization Strategy}

Training is performed in three stages. First, a warm-up stage optimizes
only the free-boundary parameter while the neural-network weights are
kept fixed. The warm-up stage consists of 2000 Adam epochs with a
learning rate of $5\times10^{-3}$. The network output layer is initially
set to zero, so that the initial structured trial solution is close to
$u(\xi)=\xi^p$.

In the second stage, all network parameters and the free-boundary
parameter are optimized jointly using Adam. This stage consists of
$10\,000$ epochs with an initial learning rate of $10^{-3}$. The learning
rate is reduced by a factor of $0.5$ when the loss does not decrease for
500 consecutive epochs. Gradient clipping with threshold $1.0$ is used
throughout the training.

Finally, L-BFGS is used for high-accuracy refinement. Two successive
refinement stages are performed using fixed collocation sets containing
4000 and 8000 points, respectively. Each L-BFGS stage uses a maximum of
500 iterations and the strong-Wolfe line search. Fixing the collocation
points during each refinement stage ensures that the optimizer acts on a
fixed discrete loss function.

\subsection{Evaluation Metrics}

The trained PINN is evaluated using the learned dead-core location and
the error in the concentration profile. The interface error is defined
as
\begin{equation}\label{eq_interface_error}
E_{x_{dz}}
=
\left|
x_{dz}^{\mathrm{PINN}}
-
x_{dz}^{\mathrm{ref}}
\right|,
\end{equation}
where $x_{dz}^{\mathrm{ref}}$ is obtained from the analytical expression
in Eq.~\eqref{eq_xdz_power}.

The global concentration error is measured using
\begin{equation}\label{eq_l2_error}
\|u_{\mathrm{PINN}}-u_{\mathrm{ref}}\|_{L^2}
=
\left(
\int_0^1
\left[
u_{\mathrm{PINN}}(x)-u_{\mathrm{ref}}(x)
\right]^2\dx
\right)^{1/2}.
\end{equation}
In the numerical experiments, the integral is evaluated using the
trapezoidal rule on a uniform grid of 1200 points.

\section{Numerical Results and Discussion}

The proposed free-boundary PINN is tested on four representative
dead-core problems with different reaction orders and Thiele moduli.
The selected cases include small and moderate reaction orders, a
reaction order close to one, and a near-critical configuration. They
therefore cover different active-region sizes and solution shapes,
including the strongly localized profile associated with a large
asymptotic exponent $p=2/(1-n)$.

For each test case, the reference dead-core interface is computed from
the analytical expression
\begin{equation}\label{eq_xdz_reference}
x_{dz}^{\mathrm{ref}}
=
1-\frac{\phi^*}{\phi},
\end{equation}
where $\phi^*$ is given by Eq.~\eqref{formula_phistar}. The accuracy of
the PINN is assessed using the absolute interface error
\[
|\Delta x_{dz}|
=
\left|
x_{dz}^{\mathrm{PINN}}-x_{dz}^{\mathrm{ref}}
\right|
\]
and the $L^2$ error of the concentration profile. The results are
summarized in Table~\ref{tab:results_table}.

\begin{table}[t]
\centering
\small
\setlength{\tabcolsep}{6pt}
\renewcommand{\arraystretch}{1.15}
\begin{tabular}{@{}cccccc@{}}
\toprule
$n$ & $\phi$ &
$x_{dz}^{\mathrm{PINN}}$ &
$x_{dz}^{\mathrm{ref}}$ &
$|\Delta x_{dz}|$ &
$L^2$ error \\
\midrule
$0.01$ & $6.0000$
& $0.760732$ & $0.760729$
& $2.28\times10^{-6}$ & $1.073\times10^{-6}$ \\

$0.50$ & $6.0000$
& $0.422652$ & $0.422650$
& $2.65\times10^{-6}$ & $1.035\times10^{-6}$ \\

$0.90$ & $29.2404$
& $0.333334$ & $0.333333$
& $8.00\times10^{-7}$ & $1.295\times10^{-5}$ \\

$0.10$ & $1.6645$
& $0.009901$ & $0.009901$
& $3.30\times10^{-7}$ & $6.806\times10^{-7}$ \\
\bottomrule
\end{tabular}
\vspace{2mm}
\caption{Accuracy of the PINN approximation for four dead-core test
cases. Interface errors are computed using the full-precision numerical values of $x_{dz}$. The interface locations displayed in the table are rounded to six decimal places.}
\label{tab:results_table}
\end{table}

Table~\ref{tab:results_table} reports the results of a representative training run for each test case, while the sensitivity to random initialization and collocation sampling is assessed separately in Table~\ref{tab:robustness}.

The results demonstrate accurate recovery of both the concentration
profiles and the free-boundary locations. For the first two cases,
$n=0.01$ and $n=0.50$ with $\phi=6$, the PINN identifies the interface
to within approximately $3\times10^{-6}$ while the $L^2$ errors are of
order $10^{-6}$. The case $n=0.90$ is more demanding because
\[
p=\frac{2}{1-n}=20.
\]
The concentration therefore remains nearly flat over most of the active
region and increases sharply close to the slab surface. Despite this
behavior, the interface error remains below $10^{-6}$ and the global
$L^2$ error remains of order $10^{-5}$. Finally, the case
$n=0.10$, $\phi=1.6645$ is close to the critical regime. Since
\[
\phi^*
=
\frac{\sqrt{2(1.1)}}{0.9}
\approx 1.6480,
\]
the corresponding dead-core length is very small,
$x_{dz}\approx0.0099$. The PINN, nevertheless, recovers this interface
with an absolute error below $10^{-6}$.

\begin{figure}[t]
    \centering
    \includegraphics[width=1\linewidth]{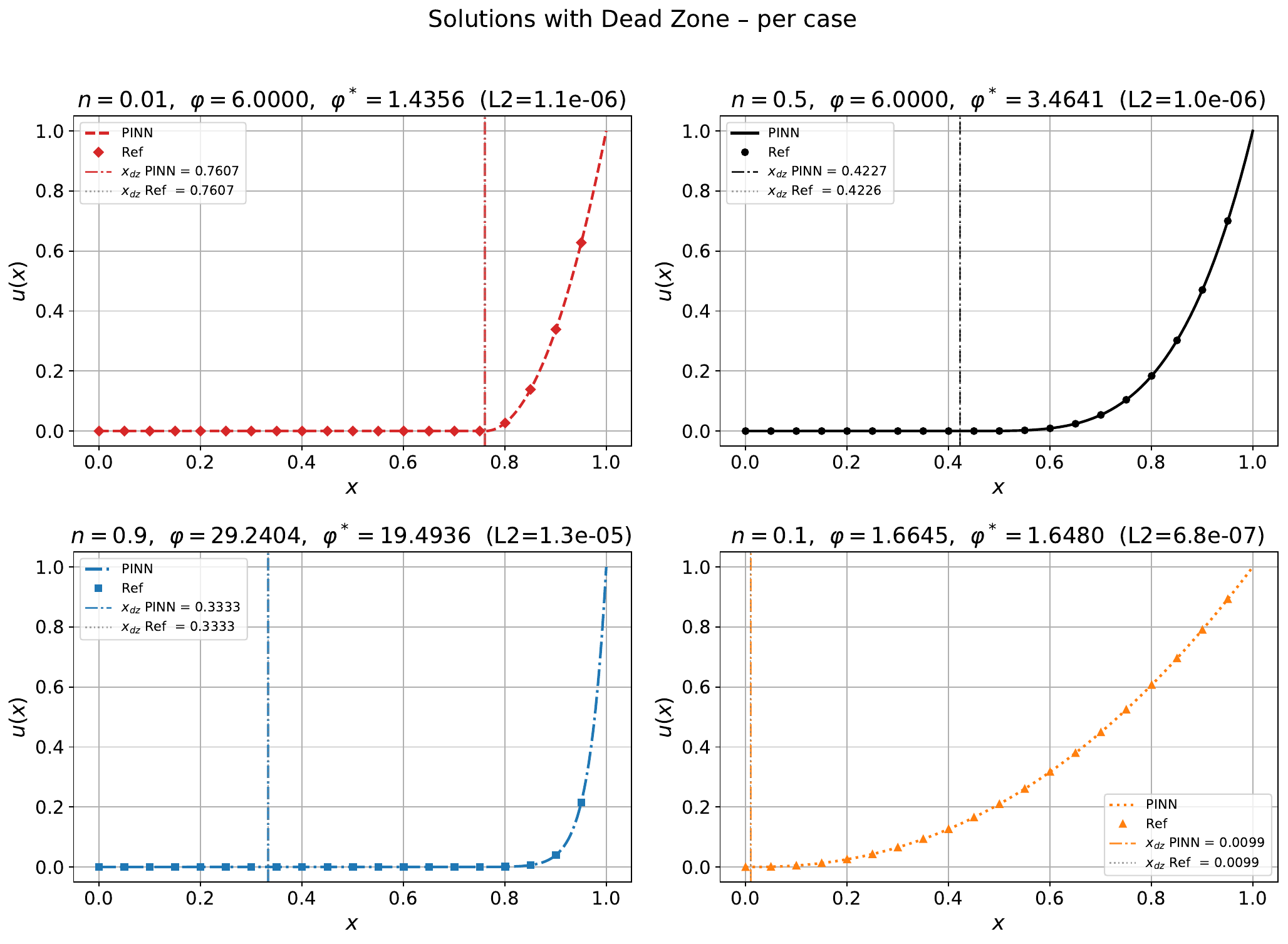}
    \caption{PINN and reference concentration profiles for
    the four dead-core test cases.}
    \label{fig:pinn_reference}
\end{figure}

Figure~\ref{fig:pinn_reference} compares the PINN predictions with the
corresponding reference solutions. The agreement is close throughout
the computational domain. In cases with pronounced dead cores, the
PINN reproduces both the zero-concentration region and the active-region
profile. For $n=0.90$, the network also captures the steep increase close
to the outer boundary, while the near-critical case exhibits a smooth
profile with only a very small dead-core region.

\begin{figure}[t]
    \centering
    \includegraphics[width=1\linewidth]{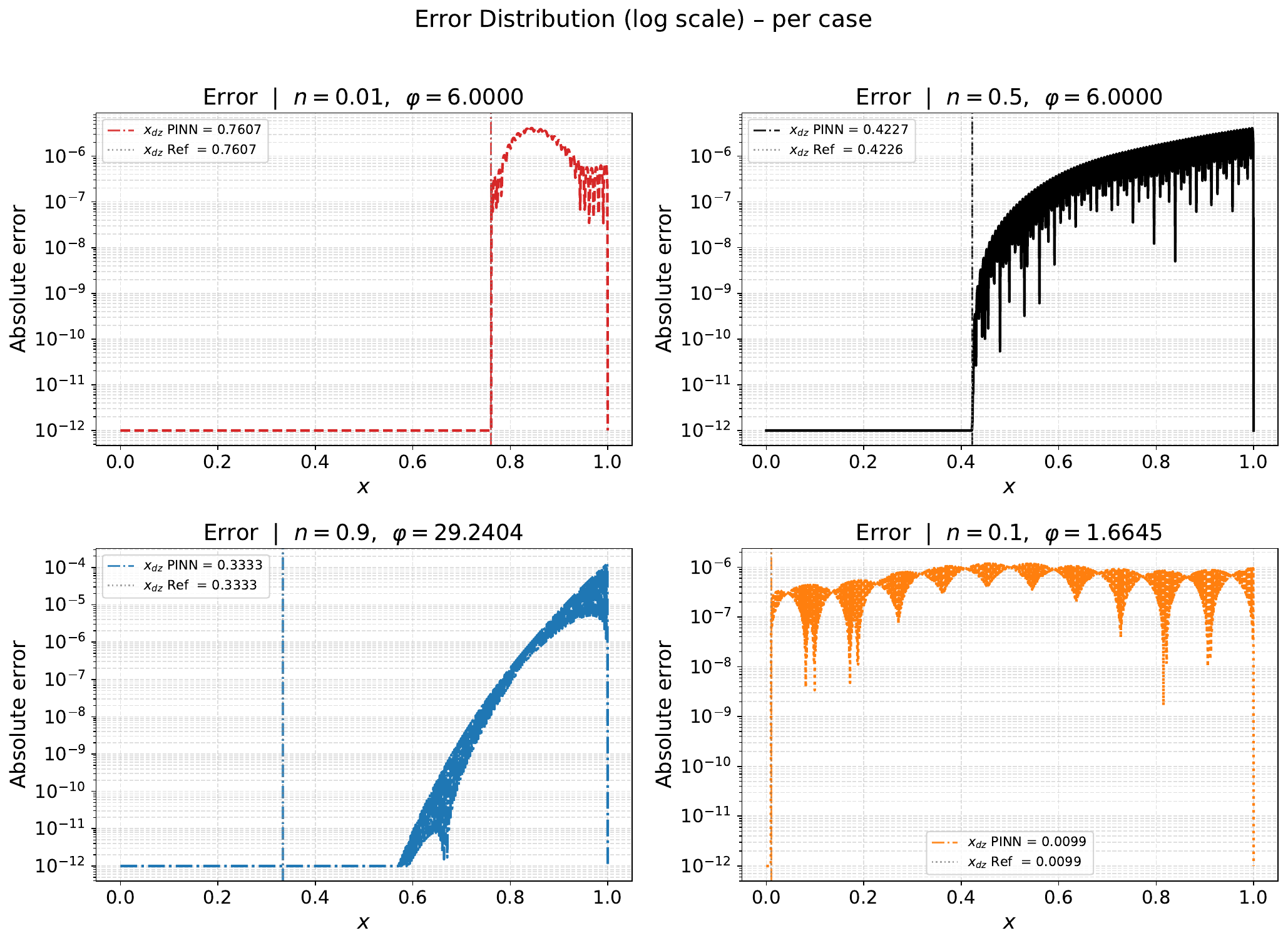}
    \caption{Pointwise absolute error
    $\left|u_{\mathrm{PINN}}(x)-u_{\mathrm{ref}}(x)\right|$
    on a logarithmic scale for the four dead-core test cases. The
    vertical lines indicate the corresponding dead-core interface
    locations.}
    \label{fig:pinn_error}
\end{figure}

To examine the spatial distribution of the approximation error, the
pointwise absolute error
\[
e(x)
=
\left|
u_{\mathrm{PINN}}(x)-u_{\mathrm{ref}}(x)
\right|
\]
is shown on a logarithmic scale in Fig.~\ref{fig:pinn_error}.

The error remains small throughout the domain. The largest deviations
occur near regions where the solution changes most rapidly. For cases
with a pronounced dead core, this occurs near the dead-core interface,
where the solution changes from zero to positive concentration. In the
large-$p$ case $n=0.90$, the error also increases near the outer boundary,
where the concentration rises sharply. The sharp local minima in the
error curves correspond to points where the PINN and reference profiles
intersect and the pointwise error becomes very small.

\subsection{Comparison with a Classical Shooting Method}

The proposed PINN is compared with a classical shooting method to provide
an independent numerical benchmark and to quantify the computational cost.
Direct shooting in the concentration variable becomes poorly conditioned
for large $p$ because
\[
u(x)\sim C(x-x_{dz})^{2/(1-n)}
\]
near the dead-core interface. For the case $n=0.90$, this exponent is
$20$. We therefore use the transformed variable
\[
y=u^{(1-n)/2},
\]
for which the leading-order interface behavior is linear. The unknown
interface location is then determined by shooting and bisection. The
shooting calculations use $\varepsilon=10^{-6}$ and the same numerical
tolerances for both test cases.

Table~\ref{tab:shooting_comparison} compares the PINN with the transformed
shooting method for a moderate case and a large-$p$ case. The interface
error is measured relative to the exact value given by
Eq.~\eqref{eq_xdz_reference}. Shooting times are averages over 20 runs
using the same hardware.

\begin{table}[t]
\centering
\small
\setlength{\tabcolsep}{6pt}
\renewcommand{\arraystretch}{1.15}
\begin{tabular}{@{}ccccr@{}}
\toprule
$n$ & $\phi$ & Method & $|\Delta x_{dz}|$ & Time (s) \\
\midrule
$0.50$ & $6.00$
& PINN & $2.65\times10^{-6}$ & $365.907$ \\
&
& Shooting & $7.45\times10^{-12}$ & $0.007$ \\
\addlinespace[1mm]
$0.90$ & $29.24$
& PINN & $8.00\times10^{-7}$ & $364.527$ \\
&
& Shooting & $5.10\times10^{-13}$ & $0.031$ \\
\bottomrule
\end{tabular}
\vspace{2mm}
\caption{Comparison of the proposed PINN with the transformed classical
shooting method for two representative dead-core cases. Shooting times
are averages over 20 runs using the same hardware.}
\label{tab:shooting_comparison}
\end{table}

As expected, shooting is substantially faster for this one-dimensional
benchmark and provides highly accurate interface locations. The PINN is
therefore not intended to replace specialized one-dimensional solvers.
Its main advantage is flexibility: the free boundary is trainable, and the
interface conditions are satisfied by construction, providing a framework
that can be extended to higher-dimensional and analytically intractable
problems.

\subsection{Robustness and Ablation Study}

To assess robustness with respect to random initialization and
collocation sampling, each of the four test cases was repeated five
times using independent random seeds. The resulting means and standard
deviations are reported in Table~\ref{tab:robustness}. No failed runs were observed. The absolute values of both the profile and interface errors remain small across all random seeds, although their relative variability is larger when the mean error approaches the numerical accuracy of the computation.

\begin{table}[t]
\centering
\small
\setlength{\tabcolsep}{6pt}
\renewcommand{\arraystretch}{1.15}
\begin{tabular}{@{}cccccc@{}}
\toprule
$n$ & $\phi$ &
\multicolumn{2}{c}{$L^2$ error} &
\multicolumn{2}{c}{$|\Delta x_{dz}|$} \\
\cmidrule(lr){3-4}\cmidrule(lr){5-6}
& & Mean & Std. & Mean & Std. \\
\midrule
$0.01$ & $6.00$
& $1.47\times10^{-6}$ & $1.12\times10^{-6}$
& $3.50\times10^{-6}$ & $2.61\times10^{-6}$ \\

$0.50$ & $6.00$
& $1.04\times10^{-6}$ & $9.17\times10^{-10}$
& $7.76\times10^{-7}$ & $8.45\times10^{-7}$ \\

$0.90$ & $29.24$
& $1.30\times10^{-5}$ & $2.33\times10^{-10}$
& $3.28\times10^{-7}$ & $2.55\times10^{-7}$ \\

$0.10$ & $1.66$
& $5.72\times10^{-7}$ & $4.80\times10^{-10}$
& $2.53\times10^{-10}$ & $1.40\times10^{-10}$ \\
\bottomrule
\end{tabular}
\vspace{2mm}
\caption{Statistical robustness over five independent random seeds.}
\label{tab:robustness}
\end{table}

A compact ablation study was also performed for two representative
cases using three independent random seeds. The results are reported in
Table~\ref{tab:ablation}. Replacing the three-component biased sampler
with uniform-only sampling increases the interface error, particularly
for the large-$p$ case $n=0.90$. This indicates that increased
collocation density near the regions of rapid variation is beneficial
for free-boundary localization. In contrast, removing the L-BFGS
refinement has only a minor effect for the two tested cases, suggesting
that the Adam stage already produces an accurate approximation for this
problem.

\begin{table}[t]
\centering
\small
\setlength{\tabcolsep}{6pt}
\renewcommand{\arraystretch}{1.15}
\begin{tabular}{@{}lccccc@{}}
\toprule
Variant & $n$ &
\multicolumn{2}{c}{$L^2$ error} &
\multicolumn{2}{c}{$|\Delta x_{dz}|$} \\
\cmidrule(lr){3-4}\cmidrule(lr){5-6}
& & Mean & Std. & Mean & Std. \\
\midrule
Full method
& $0.50$
& $1.04\times10^{-6}$ & $1.16\times10^{-9}$
& $2.58\times10^{-7}$ & $2.17\times10^{-7}$ \\
& $0.90$
& $1.30\times10^{-5}$ & $3.15\times10^{-10}$
& $2.35\times10^{-7}$ & $5.91\times10^{-8}$ \\
\midrule
Uniform sampling 
& $0.50$
& $1.04\times10^{-6}$ & $1.01\times10^{-10}$
& $7.21\times10^{-7}$ & $4.20\times10^{-7}$ \\
only
& $0.90$
& $1.30\times10^{-5}$ & $1.00\times10^{-10}$
& $2.09\times10^{-6}$ & $1.86\times10^{-6}$ \\
\midrule
No L-BFGS
& $0.50$
& $1.04\times10^{-6}$ & $1.16\times10^{-9}$
& $2.58\times10^{-7}$ & $2.17\times10^{-7}$ \\
& $0.90$
& $1.30\times10^{-5}$ & $6.13\times10^{-10}$
& $2.37\times10^{-7}$ & $2.17\times10^{-8}$ \\
\bottomrule
\end{tabular}
\vspace{2mm}
\caption{Ablation study over three independent random seeds.}
\label{tab:ablation}
\end{table}

With the present implementation, the average PINN training time is
approximately six minutes per case. All timing experiments were performed on a laptop equipped with a 13th Gen Intel Core i7-13620H processor and 32 GB of RAM. Although an NVIDIA GeForce RTX 4070 Laptop GPU was available, the reported PINN computations were performed on the CPU. This is substantially higher than
the computational cost of the classical shooting method for the present
one-dimensional benchmark. The PINN is therefore not proposed as a
replacement for specialized one-dimensional solvers. Its principal
advantage is the flexible treatment of the free boundary within a
neural-network representation. The formulation can be extended to
problems involving multidimensional catalyst geometries, inverse
identification of kinetic or transport parameters, and parametric
studies for which explicit analytical solutions are unavailable.

Overall, the numerical experiments demonstrate that the structured PINN
accurately reproduces both the concentration profile and the dead-core
interface over the tested range of reaction orders and Thiele moduli.
The near-critical and large-$p$ cases further indicate that the
formulation remains effective when the active-region geometry or
solution shape becomes challenging.

For reproducibility, the Python implementation of the PINN
framework and MATLAB codes for the classical shooting method and numerical
comparisons are available at: $\,$
\url{https://github.com/kaisartang3/dead-core-pinn}.

\section{Conclusion}

We developed a structured Physics-Informed Neural Network (PINN) for a
nonlinear diffusion--reaction problem in a catalyst slab with an unknown
dead-core interface. The proposed formulation combines the known
asymptotic behavior at the dead-core boundary with a hard-constrained
trial solution and treats the interface location as a trainable
parameter. Consequently, the concentration profile and free boundary
are approximated simultaneously, while the interface and outer boundary
conditions are satisfied by construction rather than through penalty
terms.

For power-law kinetics, the analytical critical Thiele modulus,
dead-core length, and exact concentration profile provide reference
solutions for validation. Numerical experiments covering small and
moderate reaction orders, a near-critical case, and a large-$p$ case
demonstrate accurate recovery of both the concentration profile and
dead-core location. The robustness and ablation studies further show
that the proposed sampling strategy provides reliable interface
localization, while the L-BFGS stage offers only a modest improvement
over the preceding Adam optimization for the tested cases.

Classical shooting is more computationally efficient for the present
one-dimensional benchmark, and the proposed PINN is therefore not
intended as a replacement for specialized one-dimensional solvers.
Its main advantage is the flexible treatment of the free boundary within
the neural-network representation. Future work will focus on extending
the formulation to multidimensional catalyst geometries, inverse
identification of kinetic and transport parameters, and parametric
nonlinear diffusion--reaction problems for which analytical dead-core
solutions are unavailable.

\bibliographystyle{plain}
\bibliography{refs}
\end{document}